\documentstyle[12pt]{article}
\newcommand{\ve}{\varepsilon}
\newcommand{\be}{\begin{eqnarray}}
\newcommand{\ee}{\end{eqnarray}}
\newcommand{\n}{\nonumber }
\newcommand{\oo}{\otimes}
\newcommand{\bt}{\beta}

\newcommand{\si}{\sigma}

\newcommand{\al}{\alpha}
\newcommand{\I}{I}

\begin{document}
\begin{titlepage}
\vspace{0.5in}
\begin{center}
{\Large \bf Universal ${\cal R}$--matrix for null--plane
quantized Poincar{\'e} algebra}
\end{center}
\vspace{1in}
\begin{center}
{\Large A.I. Mudrov\footnote{e-mail: aimudrov@dg2062.spb.du}}
\end{center}
\begin{center}
\vspace{1in}
 Department of Theoretical Physics,
Institute of Physics, St.Petersburg State University,
Ulyanovskaya 1,
Stary Petergof, St. Petersburg, 198904, Russia\\
\end{center}
\vspace {1.0in}
\begin{center}
  { Abstract}
\end{center}
{\footnotesize
The universal ${\cal R}$--matrix for a quantized  Poincar{\'e}
algebra ${\cal P}(3+1)$ introduced by Ballesteros {\em et al}
is evaluated. The solution is obtained as a specific case of a
formulated multidimensional generalization to the non--standard
(Jordanian) quantization of $sl(2)$.
}
\vspace {0.5in}
\begin{center}
   1997
\end{center}
\end{titlepage}
\section{Introduction}
Recently Ballesteros {\em et al}  built a quantum deformation of
the Poincar{\'e} algebra \cite{Ball1}. The quantization found was
generated
by a triangular classical $r$--matrix and, according to Drinfeld's
theory \cite{D2}, should be a twisting of $U({\cal P}(3+1))$. An
explicit form of solution  ${\cal F}$ to the twist equation
and universal matrix ${\cal R}=\tau({\cal F}^{-1}){\cal F}$ were
not given. This problem is solved in the present paper. To know
twisting element ${\cal F}$ is very important because it deforms
not only the symmetry algebra but the geometry of the space--time
as well. Twisting a universal enveloping algebra induces coherent
transformations in modules and allows to obtain important objects
automatically, for example, to construct invariant equations
and their solutions. In order to solve the problem we resort to
the theory of quantizing Lie algebras with quasi--Abelian
dual groups (semidirect product of two Abelian subgroups)
\cite{LM3,M4,M5}. In the present communication we consider first
a class of algebras along the line of that theory. That class may
be regarded as a direct generalization of the triangular or
non-standard deformation of the Borel subalgebra in $sl(2)$
\cite{Ogiv6,VladAA7}. We find general expression for twisting
elements and universal ${\cal R}$--matrices and then apply the
developed technique to the specific case of ${\cal P}(3+1)$.
\section{General consideration}
Let ${\bf L}={\bf H}\triangleleft{\bf V} $ be a Lie algebra
splitting into a semidirect sum of its two Abelian subalgebras,
with the basic elements $H_i\in {\bf H}$ and $ X_\mu\in {\bf V}$:
$$[H_j,X_\mu]=B^i_{j\mu}H_i$$
Suppose that its dual algebra ${\bf L}^*$ has the structure
 of a semidirect sum ${\bf H}^*\triangleright {\bf V}^*$ as well,
 which is  defined via a commutative set of matrices
 $(\al^i)^\mu_\nu$. To match the consistency condition upon
${\bf L}$ and ${\bf L}^*$ it is necessary to require
$$(\al^i)^\mu_\nu B^j_{k\mu} = (\al^j)^\mu_\nu B^i_{k\mu}.$$
There exists a quantization $U_\al({\bf L})$ of the universal
enveloping algebra $U({\bf L})$  with the relations on
the generators \cite{M4}
\be
[H_j,X_\mu]&=& \Bigl(\frac{e^{2\al\cdot H} - \I}
{2\al\cdot H}\Bigr)^\nu_\mu B^i_{j\nu}H_i
\label{com-rel}
\ee
and the coproduct
\be
\Delta_\al (H_i)   =  H_i \oo 1+ 1\oo H_i,\quad
\Delta_\al (X_\mu) = (e^{2\al\cdot H})^\nu_\mu \oo X_\nu+
X_\mu\oo 1.
\label{copr}
\ee
Symbol $\I$ stands for the unity matrix $\I^\mu_\nu=
\delta^\mu_\nu$, and $\al\cdot H$ means $\sum_i\al^i H_i$.
The apparent counit is $\ve(X_\mu)= \ve(H_i)=0$, and the
antipode may be readily found from the coproduct with the use
of the  defining axioms. Its explicit form is irrelevant for
our study.

We are interested only in such $U_\al({\bf L})$ which are
obtained by twisting classical universal enveloping algebras
$U({\bf L})$. To find the explicit form of element
${\cal F}\in U({\bf L})\oo U({\bf L})$ governing that process
and the universal ${\cal R}$--matrix of algebra $U_\al({\bf L})$
is the goal of our investigation. Unexpectedly, it appears easier
to start from $U_\al({\bf L})$ rather than from the classical
algebra, find a solution  $\Phi$ to the twist equation, and then
return to $U({\bf L})$ (we are going to use the group properties
of twisting  and ${\cal F}=\Phi^{-1}$ in particular \cite{D8}).

We will seek a solution to the twist equation
\be
(\Delta_\al\oo id)(\Phi)\Phi_{12}&=&( id\oo \Delta_\al)(\Phi)
\Phi_{23}
\label{tw-eq}
\ee
in the form
$$\Phi = \exp(r^{i\mu}H_i \oo X_\mu)\in
U_\al({\bf L})\oo U_\al({\bf L}).$$
The classical skew--symmetric $r$--matrix then  will be
$
 r = r^{i\mu}(X_\mu\oo H_i -  H_i \oo X_\mu),
$
and matrices $\al^i$ will be expressed through the structure
constants of ${\bf L}$ by the formula
\be
 (\al^i)^{\mu}_{\nu} = \frac{1}{2} r^{j\mu} B^i_{j\nu}.
\label{alpha}
\ee
Without any loss of generality we suppose $r$ to be
non--degenerate, since otherwise we may restrict ourselves with
the image $r({\bf L}^*)=r^t({\bf L}^*) \subset {\bf L}$, which
is a subalgebra in ${\bf L}$, and twisting a subalgebra induces
that of whole ${\bf L}$.

Calculating both sides of equation (\ref{tw-eq})
\be
(\Delta_\al\oo id) (e^{r^{i\mu} H_i\oo X_\mu}) e^{r^{i\mu}
H_i\oo X_\mu \oo 1} &=&
e^{r^{i\mu} (H_i\oo 1\oo X_\mu + 1\oo H_i\oo X_\mu) }e^{r^{i\mu}
H_i\oo X_\mu \oo 1}= \n \\
&=&e^{r^{i\mu}H_i\oo 1\oo X_\mu}
e^{r^{i\mu}H_i\oo X_\mu \oo 1 + r^{i\mu}r^{j\nu}H_i\oo
[H_j, X_\mu]\oo X_\nu }
\times\n\\
&\times & e^{r^{i\mu}1\oo H_i\oo X_\mu },\n\\
( id\oo \Delta_\al) (e^{r^{i\mu} H_i\oo X_\mu}) e^{r^{i\mu}
 1 \oo H_i\oo X_\mu} &=&
e^{r^{i\mu} (H_i\oo (e^{2\al\cdot H})^\nu_\mu \oo X_\nu+H_i\oo
X_\mu\oo 1)}
e^{r^{i\mu}   1 \oo H_i\oo X_\mu}, \n
\ee
and then comparing them with each other, taking into account
commutation properties of the generators $(H_i)$ and $(X_\mu)$
we come to condition
$$
 r^{i\mu}H_i\oo 1\oo X_\mu + r^{i\mu}r^{j\nu}H_i\oo [H_j, X_\mu]
 \oo X_\nu =
r^{i\mu}H_i\oo (e^{2\al\cdot H})^\nu_\mu \oo X_\nu,
$$
which, in its turn, yields
$
 r^{i\nu}(\delta^\mu_\nu + r^{j\mu}B_{j\nu}(H)) =
 r^{i\nu}(e^{2\al\cdot H})^\mu_\nu.
$
The latter is fulfilled provided that
$
 r^{j\mu}B_{j\nu}(H) = (e^{2\al\cdot H})^\mu_\nu- \delta^\mu_\nu
$
which holds in view of (\ref{com-rel}) and (\ref{alpha}).

Our next goal is to show that formula
$\Delta(h)=\Phi^{-1} \Delta_\al(h)\Phi$,  $h\in U_\al({\bf L})$,
defines the classical comultiplication on universal enveloping
algebra $U({\bf L})$. Due to non--degeneracy of $r$--matrix,
we can lift and drop indices: $H^\mu = r^{i\mu} H_i$,
$(\al_\mu)^\rho_\nu = \al^\rho_{\mu\nu}=
r_{i\mu} (\al^i)^\rho_{\nu}$.
Matrix $(r_{i\mu})$ is the inverse to $(r^{i\mu})$:
$r^{i\mu} r_{j\mu}=\delta^i_j$.
The relations in $U_\al({\bf L})$ take the form
$$[H^\mu,X_\nu]= (e^{2\al\cdot H}-\I)^\mu_\nu.$$
With a set of numbers $\xi^\nu$ fixed, let us introduce
entities $K^\mu$ defining them as
\be
K^\mu&=&\xi^\nu(\I-e^{-2\al\cdot H})^\mu_\nu
\label{cl-bas}
\ee
and evaluate commutation relations between $K^\mu$ and $X_\nu$:
$$[K^\mu,X_\nu]=\xi^\bt (e^{-2\al\cdot H})^\rho_\bt
(2\al^\mu_{\rho\si})(e^{2\al\cdot H}-\I)^\si_\nu.$$
>From commutativity of matrices $\al_\mu$ and condition
$\al^\rho_{\mu\nu}=\al^\rho_{\nu\mu}$ following from
the classical Yang--Baxter equation we find
$\al^\mu_{\rho\si}(e^{2\al\cdot H}-\I)^\si_\nu =
\al^\mu_{\nu\si}(e^{2\al\cdot H}-\I)^\si_\rho.
$
This is verified
by a simple induction over powers of matrix $(\al\cdot H)$.
Finally, we have $$[K^\mu,X_\nu]=2\al^\mu_{\si\nu} K^\si,$$
i.e. the classical commutation relations. The coproduct on the
new generators is
$$\Delta_\al(K^\mu) = K^\mu\oo1+
(e^{-2\al\cdot H})^\mu_\nu\oo K^\nu .$$
With the use of this formula we calculate twisted coproduct which
comes out to be
\be
\Delta(K) &=&
e^{-H\oo X}\Delta_\al(K)e^{H\oo X}=\n\\
&=& K\oo1+ (e^{2\al\cdot H}e^{-2\al\cdot H})\oo K
 = K\oo1+1\oo K, \n\\
\Delta(X) &=&e^{-H\oo X}\Delta_\al(X)e^{H\oo X}=
  e^{2\al\cdot H}\oo X + e^{-H\oo X} (X\oo 1) e^{H\oo X}=\n\\
  &=&e^{2\al\cdot H}\oo X + X\oo 1 - (e^{2\al\cdot H}-\I)\oo X=
  1\oo X + X\oo 1.
\ee
We might consider our goal achieved were we sure that the number
of independent generators $K^\mu$ be the same as the
dimensionality of space ${\bf H}$. That is not the case in
general,  and it is determined by a particular choice of
$\xi^\mu$.  In the  classical limit we have
$$
K^\mu=\xi^\nu(2\al\cdot H)^\mu_\nu = [H^\mu,\xi^\nu X_\nu].
$$
Thus, while $\xi^\mu$ takes all possible values lineal
Span($K^\mu$) fills up the subspace
${\bf H}'=[{\bf H},{\bf V}]\subset {\bf H}$.
If that subspace coincides with whole ${\bf H}$  we
may state that, indeed, twisting $U_\al({\bf L})$
with element $\Phi$ results
in  $U({\bf L})$. Let us show that  dim(${\bf H}'$) $<$
dim(${\bf H}$) if and only if subalgebra ${\bf V}$ and the
center of ${\bf L}$ have a nontrivial intersection. Indeed,
because of non--degeneracy  of the classical  $r$--matrix there
exists an isomorphism between linear spaces  ${\bf V}^*$ and
${\bf H}$ (basic elements  $H^\mu $ and $X_\mu$  turn out to be
mutually dual). Subspace ${\bf H}'$ is less than ${\bf H}$ if and
only if there is an element $X_0=\xi^\mu_0 X_\mu \in {\bf H}^*$,
orthogonal to entire ${\bf H}'$:
$$0=(X_0,[H^\mu,X_\nu])=2 \xi_0^\si \al^\mu_{\si\nu}.$$
Due to the lower indices symmetry of tensor $\al^\mu_{\si\nu}$,
the latter expression is nothing else than the matrix of the
adjoint representation ad$(X_0)$ restricted to subspace ${\bf V}$.
Let us summarize the results of our study.
\newtheorem{Rmatrix}{Theorem}
\begin{Rmatrix}
Element $\Phi = \exp(r^{i\mu}H_i \oo X_\mu)\in
U_\al({\bf L})\oo U_\al({\bf L}) $
is a solution to the twist equation.
Twisting $U_\al({\bf L})$ with the help of $\Phi$ gives classical
universal enveloping algebra $U({\bf L})$ unless ${\bf V}$ contains
central elements. The universal ${\cal R}$--matrix for
$U_\al({\bf L})$ is
\be
 {\cal R} &=& \exp(r^{i\mu}X_\mu\oo H_i)\exp(-r^{i\mu}
 H_i \oo X_\mu)
 \label{R-mat}
\ee
\end{Rmatrix}
Expression (\ref{R-mat}) is an apparent generalization of
the ${\cal R}$--matrix for $U_h(sl(2))$ in the form by O.
Ogievetsky.

It becomes clear from the above study that the algebras
$U_\al({\bf L})$ covered by the theorem
are completely specified by the set of
commutative matrices $\al_\mu$, satisfying
requirement $(\al_\mu)^\si_\nu=(\al_\nu)^\si_\mu$.
Such matrices define an associative commutative multiplication
$X_\mu\circ X_\mu \equiv \al^\si_{\mu\nu} X_\si$
on the subspace ${\bf V}$ and  {\em vice versa}. For an
example of this kind let us take
$(\al_\mu)^\si_\nu= \delta^\si_{\mu+\nu}$,
where  $\mu,\nu,\si=0,\ldots,n $.
Thus introduced, $\al_\mu$ form an Abelian matrix ring and
define a semidirect sum ${\bf L}={\bf H}\triangleleft{\bf V}$
with ${\bf H}\sim {\bf V}^* $. It complies with the condition
of the theorem since matrix $\al_0$, the unity of the
ring, has the maximum rank equal
to $n+1$, and we may set $\xi^\mu=\delta^0_\mu$ in transformation
(\ref{cl-bas}).

\section{Universal ${\cal R}$--matrix for the
quantum Poincar{\'e} algebra} Let us apply the developed
technique to the quantum universal enveloping Poincar{\'e}
algebra. The deformation is generated by twisting of the
subalgebra
${\bf L}$=Span($E_1$, $E_2$, $P_1$, $P_2$, $P_+$, $K_3$),
in terms of \cite{Ball1}. Notations $H^i=-z E_i$, $H^3=-z P_+$,
$Y_i=2P_i$, $Y_3=-2K_3$, $i=1,2$, having been  introduced, the
coproduct in $U_\al({\bf L})$ reads:
\be
\Delta_\al (H^\mu)   &=&  H^\mu \oo 1+ 1\oo H^\mu,\n\\
\Delta_\al (Y_1) &=& e^{H^3 } \oo Y_1+ Y_1\oo e^{-H^3 },\n\\
\Delta_\al (Y_2) &=& e^{H^3 } \oo Y_2+ Y_2\oo e^{-H^3 },\n\\
\Delta_\al (Y_3) &=& e^{H^3 } \oo Y_3+ Y_3\oo e^{-H^3 }+\n\\
               &+& e^{H^3 }H^1 \oo Y_1- Y_1\oo H^1 e^{-H^3 }+\n\\
               &+& e^{H^3 }H^2 \oo Y_2- Y_2\oo H^2 e^{-H^3 }.\n
\ee
Non--vanishing commutators are
$$
[H^i,Y_i]  = 2\sinh (H^3),\quad [H^i,Y_3]  = 2\cosh (H^3)H^i,
\quad [H^3,Y_3]  = 2\sinh(H^3),\quad i=1,2.
$$
Correspondence with the notations of the previous paragraph
is achieved through the transformation
$X_\mu = Y_\nu (e^{\al\cdot H})_\mu^\nu$, the matrices
$\al$ being
$$
\al_1=\left(
\begin{array}{lll}
0 & 0 & 1\\
0 & 0 & 0\\
0 & 0 & 0\\
\end{array}
\right),\quad
\al_2=\left(
\begin{array}{lll}
0 & 0 & 0\\
0 & 0 & 1\\
0 & 0 & 0\\
\end{array}
\right),\quad
\al_3=\left(
\begin{array}{lll}
1 & 0 & 0\\
0 & 1 & 0\\
0 & 0 & 1\\
\end{array}
\right).
$$
Explicitly this results in the following change of variables:
$$X_1=Y_1e^{H^3},\quad X_2=Y_2e^{H^3},\quad
X_3=(Y_1 H^1 + Y_2 H^2 + Y_3)e^{H^3}.$$
Transition to the classical generators is completed by
transformation  (\ref{cl-bas}), where we may assume
$\xi^1=\xi^2=0$, $\xi^3={1\over 2}$, for  matrix $\al_3$ has the
maximum rank 3. $$K^1 = H^1 e^{-2H^3},\quad K^2 = H^2 e^{-2H^3},\quad
K^3 = {1\over2}(1-e^{-2H^3}).$$
Elements $K^\mu$ and $X_\nu$ obey the ordinary, non--quantum,
commutation relations of  $U({\bf L})$:
$$
[K^i,X_i]  = 2 K^3,\quad [K^i,X_3]  = 2 K^i,
\quad [K^3,X_3]  = 2 K^3,\quad i=1,2.
$$
Quantization $U({\bf L})\to U_\al({\bf L})$ is controlled by the
twisting element
$${\cal F}=\exp\Bigl(
\frac{K^1}{2K^3-1}\oo X_1
+\frac{K^2}{2K^3-1}\oo X_2
+\frac{1}{2}\ln(1-2 K^3)\oo X_3 \Bigr)= \exp({-H^\mu\oo X_\mu}),$$
and the quantum universal ${\cal R}$--matrix of the algebra
$U_\al({\bf L})$ is given by
$$
 {\cal R} = \exp(X_\mu\oo H^\mu)\exp(-H^\mu \oo X_\mu).
$$
\section{Conclusion}
The present investigation continues the series of works
\cite{LM3,M4,M5} devoted to a method of constructing quantum Lie
algebras with the use of  a classical object, the dual group.
Based on the duality principle \cite{STS9} viewing a quantum
universal Lie algebra as a set of non--commutative functions on
the dual group, that method reduces the quantizaton problem to
finding a deformed Lie biideal consistent with a given coproduct
\cite{M4}. Because of complicated  structure of a generic Lie group,
that program is feasible, however, for simple classes or in separate
particular cases. So is the set of quasi--Abelian groups which are
decomposed into a semidirect composition of two Abelian subgroups.
Quantization theory for such Lie algebras was developed in
\cite{M4}. Simple as dual groups of that type may seem, they occur
rather frequently, especially in low dimensions
\cite{Ball1,Ball10,Ball11,Ball12,Ball13}, and corresponding quantum
algebras possess very diverse and nontrivial properties. So, that
class includes two non--isomorphic deformations of $U(sl(2))$. A
generalization to the standard quantization was studied in detail
in our work \cite{M5}. Here we have considered a generalization to
the non--standard quantization of $sl(2)$ or, to be more exact, its
Borel subalgebra. The universal ${\cal R}$--matrix for the Jordanian
$U_h(sl(2))$ has been found by O. Ogievetsky and A. A. Vladimirov
\cite{Ogiv6,VladAA7,VladAA14}. In the present paper we have obtained
explicit expression for ${\cal F}$ and  ${\cal R}$ for the
multi--dimensional generalization of the Borel subalgebra in $sl(2)$.
The technique developed made it possible to built the universal
${\cal R}$--matrix for the null--plane quantizaton of the
Poincar{\'e} algebra ${\cal P}(3+1)$ found in \cite{Ball1}.

\vspace{0.3cm}
\noindent
{\Large\bf  Acknowledgement}
\vspace{0.3cm}

\noindent
We are grateful to A. A. Vladimirov who informed us about the
work by O. Ogievetsky  \cite{Ogiv6}.

\end{document}